\documentclass[12pt]{amsart}
\usepackage{amssymb,latexsym,eufrak,amsmath,amscd, graphics}
\usepackage{xspace,xcolor}
\usepackage[breaklinks,colorlinks=true,pagebackref,hyperindex]{hyperref}
\usepackage[alphabetic]{amsrefs} 
\usepackage{mathrsfs}

\usepackage{amsfonts}

\usepackage{tikz-cd} 
\usepackage{enumerate}
\usepackage{dsfont} 
\usepackage{mathtools} 

\usepackage{stmaryrd} 

\usetikzlibrary{calc, intersections} 

\tikzset{symbol/.style={draw=none,every to/.append style={edge node={node [sloped, allow upside down, auto=false]{$#1$}}}}}

\newcommand{\into}{\hookrightarrow}

\DeclareMathOperator{\Sp}{Sp}

\DeclareMathOperator{\colim}{\operatorname{colim}}

\DeclareMathOperator{\lcm}{\operatorname{lcm}}

\DeclareMathOperator{\lct}{\operatorname{lct}}

\DeclareMathOperator{\Der}{\operatorname{Der}}

\DeclareMathOperator{\Ann}{\operatorname{Ann}}

\DeclareMathOperator{\Gr}{\operatorname{Gr}}

\DeclareMathOperator{\Spec}{\operatorname{Spec}}

\newtheorem{lemma}{Lemma}[section]
\newtheorem{theorem}[lemma]{Theorem}

\theoremstyle{definition}

\newtheorem{example}[lemma]{Example}

\theoremstyle{remark}
\newtheorem*{remark*}{Remark}
\newtheorem*{note*}{Note}

\let\oldtocsection=\tocsection
\let\oldtocsubsection=\tocsubsection
\let\oldtocsubsubsection=\tocsubsubsection
\renewcommand{\tocsection}[2]{\hspace{0em}\oldtocsection{#1}{#2}}
\renewcommand{\tocsubsection}[2]{\hspace{2em}\oldtocsubsection{#1}{#2}}
\renewcommand{\tocsubsubsection}[2]{\hspace{3em}\oldtocsubsubsection{#1}{#2}}

\begin{document}

\title{Multiplicative Thom-Sebastiani for Bernstein-Sato polynomials}

\begin{abstract}
We show that if $f\in \mathcal{O}_X(X)$ and $g\in \mathcal{O}_Y(Y)$ are nonzero regular functions on smooth complex algebraic varieties $X$ and $Y$, then the Bernstein-Sato polynomial of the product function $fg \in \mathcal{O}_{X\times Y}(X \times Y)$ is given by $b_{fg}(s)=b_f(s)b_g(s)$, answering a question of Budur and Popa.
\end{abstract}

\thanks{The author was partially supported by a Rackham Merit Fellowship.}


\author[Jong]{Jonghyun Lee}
\address{Department of Mathematics, University of Michigan,
Ann Arbor, MI 48109, USA}
\email{{nuyhgnoj@umich.edu}}

\maketitle

\markboth{Jonghyun Lee}{Multiplicative Thom-Sebastiani for Bernstein-Sato polynomials}


\section{Introduction}

Let $f\in \mathcal{O}_X(X)$ be a nonzero regular function on a smooth complex algebraic variety $X$. An important measure of singularity of $f$ is its Bernstein-Sato polynomial $b_f(s) \in \mathbb{C}[s]$, which is the monic polynomial of minimal degree satisfying the relation
$$
b_f(s) \cdot f^s \in \mathcal{D}_X[s]\cdot f^{s+1}.
$$
Here $s$ is an indeterminate and $\mathcal{D}_X$ is the sheaf of differential operators on $X$, acting on the formal power $f^s$ in the expected way. Independently discovered by Bernstein \cite{Ber72} and Sato, with a general existence proof in \cite{Bjo93} (in the analytic setting), the Bernstein-Sato polynomial is connected to several other invariants of singularities. For example, the roots of $b_f(s)$ determine the eigenvalues of the monodromy action on the cohomology of the Milnor fiber by a result of Malgrange \cite{Mal83}, and the negative of the largest root of $b_f(s)$ is equal to the log canonical threshold $\lct(f)$ of $f$ by a result of Lichtin and Koll\'{a}r \cite[Theorem 10.6]{Kol97}. 



The Bernstein-Sato polynomial is a subtle invariant of $f$, and even seemingly basic properties resist immediate proofs. For example, while it is easy to show that $b_f(s)=s+1$ if $f$ cuts out a smooth hypersurface, the converse, proven first in \cite[Proposition 2-6]{BM96} and later shown to follow from properties of Hodge ideals in \cite[Remark 6.4]{MP20}, has yet to admit an elementary proof. 

In this paper we show the following seemingly basic property, a so-called multiplicative Thom-Sebastiani formula for Bernstein-Sato polynomials.

\begin{theorem}\label{main thm}
Let $X,Y$ be smooth complex algebraic varieties, $f\in \mathcal{O}_X(X),g\in \mathcal{O}_Y(Y)$ nonzero regular functions, and consider the product function $fg=f\otimes g \in \mathcal{O}_{X \times Y}(X \times Y)$. Then
$$
b_{fg}(s) = b_f(s) b_g(s).
$$
\end{theorem}

The division $b_{fg}(s) \ | \ b_f(s) b_g(s)$ is immediate from the definition, but the reverse direction relies on a linear algebra result in Lemma \ref{main lemma} below. This answers a question of Budur in \cite[\S2.11 Questions]{Bud12} and of Popa in \cite[Remark 2.2.10 (5)]{Pop21}. The strategy of proof is to filter certain $D$-modules by finitely generated $\mathbb{C}[s]$-modules and then use the Smith normal form to explicitly work out the minimal polynomials.


We mention that another proof of Theorem \ref{main thm} has also been obtained independently by Shi and Zuo in \cite[Theorem 3.2]{SZ24}. Both our proofs begin in a similar fashion by filtering the same $D$-modules by $\mathbb{C}[s]$-submodules. The main difference in our methods essentially comes down to how we compute minimal polynomials: while we filter by finite free modules and then use explicit bases given by the Smith normal form together with flatness of certain quotient modules, Shi and Zuo filter by possibly non-finitely generated free modules and then use a general result on intersections of tensor products of free modules together with the associated graded of regular local rings.

We note that the additive Thom-Sebastiani formula for Bernstein-Sato polynomials (in the case one of the functions is quasihomogeneous) was worked out in \cite[Proposition 0.7, Theorem 0.8]{Saito_microlocal}.



\subsection{Outline of the paper}

In \S\ref{setup} we set up the notation and prove Theorem \ref{main thm}, which is an application of a linear algebra lemma shown in \S\ref{main lemma sec}. In the last section \S\ref{additional remarks}, we make additional remarks on Bernstein-Sato ideals, Milnor fiber monodromy, $V$-filtrations, and generalized $b$-functions.

\subsection{Acknowledgements}

I'd like to thank Bradley Dirks for helpful discussions on Bernstein-Sato polynomials and $V$-filtration, as well as for providing comments on this paper. I am indebted to Uli Walther for bringing to my attention the reference \cite{Tap08} and to Nero Budur for bringing to my attention his survey \cite{Bud12}. I am grateful to my advisor Mircea Musta\c{t}\u{a} for our weekly meetings and his careful readings and comments on this paper. I'd also like to thank the referee for several helpful comments.

\section{The setup and proof of theorem}\label{setup}

To prove Theorem \ref{main thm}, we may assume $X,Y$ are affine because $b_f(s) = \lcm_i b_{f|_{U_i}}(s)$ for an affine open cover $X=\cup_i U_i$ of $X$. So let $X=\Spec R$ and $Y=\Spec S$ be smooth affine complex varieties, with $f \in R, g\in S$. Consider the product $fg=f\otimes g\in R \otimes_{\mathbb{C}}S$. 


For basic facts on rings of differential operators and $\mathcal{D}$-modules, see \cite{HTT08}. Let $D_X:=\Gamma(X, \mathcal{D}_X)$ be the ring of $\mathbb{C}$-linear differential operators of $R$, and similarly define $D_Y,D_{X\times Y}$. Consider the following module
$$
B_{f,+}:= R_f[s]f^s,
$$
where $s$ is an indeterminate and $B_{f,+}$ is a free module over the polynomial ring $R_f[s]$ with free generator the formal symbol $f^s$. We denote $f^{s+m}=f^m f^s$ for $m\in \mathbb{Z}$. $B_{f,+}$ has a left module action by the ring $D_X[s]$ of differential operators on $X$ adjoined the formal variable $s$, where $s$ commutes with $D_X$. The $s$-action on $B_{f,+}$ comes from the $R_f[s]$-module structure, and the $D_X$-action is given by
$$
\theta \cdot u s^j f^s = \theta(u) s^jf^s + \frac{\theta(f)}{f} us^{j+1} f^s,
$$
for $\theta\in \Der_{\mathbb{C}}(R)$, $u\in R_f$, and $j \in \mathbb{Z}_{\geq 0}$. We can express the Bernstein-Sato polynomial of $f$ as follows:
$$
\Ann_{\mathbb{C}[s]} (D_X[s] f^s/D_X[s] f^{s+1}) = (b_f(s)).
$$
Form similar modules for $g$ and $fg$. Using
$$
(R\otimes_{\mathbb{C}} S)_{fg}=R_f\otimes_R R\otimes_{\mathbb{C}} S \otimes_SS_g = R_f \otimes_{\mathbb{C}} S_g,
$$
we obtain the following isomorphism linear over $D_{X\times Y}[s] =D_X[s] \otimes_{{\mathbb{C}}[s]} D_Y[s]$
\begin{align}
B_{fg,+}=(R\otimes_{\mathbb{C}} S)_{fg}[s](fg)^s &\simeq R_f[s]f^s \otimes_{{\mathbb{C}}[s]} S_g[s]g^s = B_{f,+} \otimes_{{\mathbb{C}}[s]} B_{g,+} \label{tensor} \\
(u\otimes v) s^j (fg)^s &\leftrightarrow us^jf^s \otimes v g^s = uf^s \otimes vs^j g^s \nonumber,
\end{align}
for $u\in R_f, v\in S_g, j\in \mathbb{Z}_{\geq 0}$. The following commutative diagram
\[
\begin{tikzcd}
D_{X\times Y}[s] \arrow[r, symbol=\simeq]\arrow[d, "(fg)^s"] & D_X[s] \otimes_{{\mathbb{C}}[s]} D_Y[s] \arrow[d, "f^s\otimes g^s"] \\
B_{fg,+} \arrow[r, symbol=\simeq] & B_{f,+} \otimes_{\mathbb{C}[s]} B_{g,+},
\end{tikzcd}
\]
induces an isomorphism
$$
D_{X\times Y}[s](fg)^s \simeq D_X[s]f^s \otimes_{{\mathbb{C}}[s]} D_Y[s]g^s.
$$
Similarly, we have 
$$
D_{X\times Y}[s](fg)^{s+1} \simeq D_X[s]f^{s+1} \otimes_{{\mathbb{C}}[s]} D_Y[s]g^{s+1}.
$$

To apply Lemma \ref{main lemma} below, note that $R_f[s]f^s \simeq R_f[s]$ as a ${\mathbb{C}}[s]$-module, hence $R_f[s]f^s$ is ${\mathbb{C}}[s]$-torsion-free. Moreover, observe that $D_X[s]f^s/{\mathbb{C}}[s]f^s$ is ${\mathbb{C}}[s]$-torsion-free because 
$$
D_X[s]f^s/{\mathbb{C}}[s]f^s \subset R_f[s]f^s/{\mathbb{C}}[s]f^s \simeq_{{\mathbb{C}}[s]} R_f[s]/{\mathbb{C}}[s],
$$
where $R_f[s]/{\mathbb{C}}[s]$ is ${\mathbb{C}}[s]$-torsion-free because we can split the inclusion $\mathbb{C}[s] \into R_f[s]$ by evaluating at a closed point $x \in \Spec R$ where $f$ does not vanish (in fact, $R_f[s]/{\mathbb{C}}[s]$ is ${\mathbb{C}}[s]$-free because $1\in R_f$ being a part of a ${\mathbb{C}}$-basis of $R_f$ implies that $1$ is also a part of a ${\mathbb{C}}[s]$-basis of $R_f[s]$).

\noindent{\it Proof of Theorem \ref{main thm}}. Take $k=A_1=A_2 :={\mathbb{C}}[s],P_1:=D_X[s]f^s, Q_1:=D_X[s]f^{s+1},P_2:=D_Y[s]g^s, Q_2:=D_Y[s]g^{s+1}, \delta_1:=f^s,\delta_2:=g^s,b_1:=b_f(s),b_2:=b_g(s)$ in Lemma \ref{main lemma} below, and conclude.  \qed

\section{Main lemma}\label{main lemma sec}

Here is the setup for the main linear algebra lemma. Fix $r\geq 1$ and a principal ideal domain $k$. Let $A_1,\dots,A_r$ be principal ideal domains flat over $k$ such that the tensor product $A=A_1\otimes_k\cdots \otimes_k A_r$ is a unique factorization domain. Let $Q_i \into P_i$ be an inclusion of flat $A_i$-modules for $i=1,\dots,r$. Let $\delta_i \in P_i$ be nonzero elements such that 
\begin{enumerate}[(1)]
\item their images $\overline{\delta_i}$ in $P_i/Q_i$ are $A_i$-torsion and
\item the quotients $P_i/A_i\delta_i$ are $A_i$-flat.
\end{enumerate}
Let $P=P_1\otimes_k\cdots\otimes_k P_r$,  $Q=Q_1\otimes_k \cdots \otimes_k Q_r$, and $\delta = \prod_{i=1}^r\delta_i \in P$, whose image in $P/Q$ we denote by $\overline{\delta}$. 

For the following lemma, recall that being torsion-free and being flat are equivalent over a principal ideal domain.

\begin{lemma}\label{main lemma}
If $\Ann_{A_i}(\overline{\delta_i})=(b_i)$ for $b_i\in A_i$ for $1\leq i\leq r$, and $b=b_1\otimes \cdots \otimes b_r \in A$, then 
$$
\Ann_A(\overline{\delta}) = (b).
$$
\end{lemma}
\noindent{\it Proof}. Write $P_i=\cup_{j_i} P_{ij_i}$ as a filtered union of its finitely generated submodules containing $\delta_i$. Let $Q_{ij_i} = Q_i \cap P_{ij_i}$, so that $Q_i=\cup_j Q_{ij_i}$ and we get the filtered unions
$$
P = \bigcup_{j_1,\dots,j_r}\otimes_i P_{ij_i}, \ \ \ Q =  \bigcup_{j_1,\dots,j_r}\otimes_i Q_{ij_i},  \ \ \ P_i/Q_i = \bigcup_{j_i} P_{ij_i}/Q_{ij_i}.
$$
Note that the inclusion $Q \into P$ is a filtered union of the inclusions $\otimes_i Q_{ij_i} \into \otimes_i P_{ij_i}$, so because colimits commutes with cokernels, we have
$$
P/Q = \colim_{j_1,\dots,j_r} (\otimes_iP_{ij_i}/\otimes_iQ_{ij_i}).
$$
By construction, we can consider the image $\overline{\delta}_{j_1,\dots,j_r}$ of $\delta$ in $\otimes_iP_{ij_i}/\otimes_iQ_{ij_i}$. To show $\Ann_A(\overline{\delta}) = (b)$, it suffices to show that
$$
\Ann_A (\overline{\delta}_{j_1,\dots,j_r}) = (b)
$$
for every $j_1,\dots,j_r$. Therefore we may assume each $P_i$ is a finitely generated $A_i$-module (here we use that $k$ is a principal ideal domain which implies that the $P_{ij_i}$'s are flat over $k$). By Smith normal form, we may write the inclusion $Q_i \into P_i$ as $A^{n_i}_i \into A^{m_i}_i$ defined on standard bases by $e_p\mapsto a_{ip}e_p$. Consider the commutative diagram with exact rows and injective columns
\[
\begin{tikzcd}
0 \arrow[r] & A_i \arrow[r, "b_i"]\arrow[d,hook] & A_i \arrow[r]\arrow[d,hook, "\delta_i"] & A_i/b_i \arrow[r]\arrow[d, hook] & 0 \\
0 \arrow[r]& Q_i \arrow[r]\arrow[d,equals] & P_i \arrow[r]\arrow[d, equals] & P_i/Q_i \arrow[r]\arrow[d, equals] & 0 \\
& A_i^{n_i} & A_i^{m_i} & \bigoplus_{p=1}^{n_i} A_i/a_{ip} \oplus A_i^{m_i - n_i}. &
\end{tikzcd}
\]
Note that $A/b_i =A\delta_i \subset (P_i/Q_i)_{tors} \subset \oplus_{p=1}^{n_i} A_i/a_{ip}$, hence we may write
$$
\delta_i = (c_{i1},\dots, c_{in_i},0,\dots,0) \in A_i^{n_i}=P_i.
$$
Observe that
$$
b_i = \lcm_p (  \lcm(c_{ip}, a_{ip}) / c_{ip}  )  ,
$$
where we run through those $p$ such that $c_{ip}$  is nonzero, in which case $a_{ip}$ is nonzero because $\delta_i$ is torsion in $P_i/Q_i$. Note that since $P_i/A_i\delta_i$ is $A_i$-flat, $\delta_i$ is part of a basis of $P_i$, hence
$$
\gcd_p(c_{ip})=1.
$$
Now consider the following commutative diagram with exact rows and injective columns
\[
\begin{tikzcd}
& & A \arrow[r]\arrow[d,hook, "\delta"] & A/I \arrow[r]\arrow[d, hook] & 0 \\
0 \arrow[r]& Q \arrow[r]\arrow[d,equals] & P \arrow[r]\arrow[d, equals] & P/Q \arrow[r]\arrow[d, equals] & 0 \\
& A^{n} & A^{m} & \bigoplus_{i=1}^r\bigoplus_{p_i=1}^{n_i} A/a_{1p_1}\cdots a_{rp_r} \oplus A^{m - n} &
\end{tikzcd},
\]
where $n=n_1\cdots n_r,m=m_1\cdots m_r$, and $I=\Ann_A(\overline{\delta})$. Note that
$$
\delta=(c_{1p_1}\cdots c_{rp_r})_{i,p_i} \in A^{m}.
$$
Therefore, because $A$ is a unique factorization domain, $I$ is a principal ideal generated by
$$
l := \lcm_{p_1,\dots,p_r} (  \lcm(c_{1p_1}\cdots c_{rp_r}, a_{1p_1}\cdots a_{rp_r}) / c_{1p_1}\cdots c_{rp_r})  .
$$
To finish the proof, we need to show that $(l)=(b)\subset A$. By taking prime factorizations in the unique factorization domain $A$ and comparing the exponents of a fixed prime, it suffices to prove the following statement: 

For each $i=1,\dots,r$, let $\gamma_{ip_i},\alpha_{ip_i}\geq 0$ be finitely many nonnegative integers indexed by $p_i$ such that $\min_{p_i}(\gamma_{ip_i})=0$. Then
\begin{equation}\label{first}
\max_{p_1,\dots,p_r} (  \max(\gamma_{1p_1}+\cdots +\gamma_{rp_r}, \alpha_{1p_1}+\cdots +\alpha_{rp_r}) - \gamma_{1p_1}-\cdots -\gamma_{rp_r})  
\end{equation}
is equal to
\begin{equation}\label{second}
\sum_{i=1}^r \max_{p_i}( \max(\gamma_{ip_i}, \alpha_{ip_i}) - \gamma_{ip_i} )
\end{equation}
$$=\max_{p_1,\dots,p_r}( \max(\gamma_{1p_1},\alpha_{1p_1})+\cdots + \max(\gamma_{rp_r},\alpha_{rp_r}) - \gamma_{1p_1}-\cdots -\gamma_{rp_r}).
$$
It is immediate that (\ref{second}) is $\geq$ (\ref{first}), so we prove the reverse inequality. We need to show that for every fixed $p_1,\dots,p_r$, there exists $q_1,\dots,q_r$ such that
$$
\max(\gamma_{1p_1},\alpha_{1p_1})+\cdots + \max(\gamma_{rp_r},\alpha_{rp_r}) - \gamma_{1p_1}-\cdots -\gamma_{rp_r}
$$
$$
 \leq
 \max(\gamma_{1q_1}+\cdots +\gamma_{rq_r}, \alpha_{1q_1}+\cdots +\alpha_{rq_r}) - \gamma_{1q_1}-\cdots -\gamma_{rq_r}.
$$
If $\gamma_{ip_i}\leq \alpha_{ip_i}$, set $q_i=p_i$, and if  $\gamma_{ip_i} > \alpha_{ip_i}$, choose $q_i$ such that $\gamma_{iq_i}=0$, and we are done. \qed

\section{Additional remarks}\label{additional remarks}

\subsection{Bernstein-Sato ideals}

Let $X$ be a smooth complex algebraic variety, and let $F=(f_1,\dots,f_r): X \to \mathbb{A}^r$ be a nonconstant morphism. The Bernstein-Sato ideal $B_F \subset \mathbb{C}[s_1,\dots,s_r]$ of $F$, introduced and proven to be nonzero by Sabbah in \cite{Sab87b}, is defined as the $\mathbb{C}[s_1,\dots,s_r]$-module annihilator of the following module
$$
\mathcal{D}_X[s_1,\dots,s_r]\cdot f_1^{s_1}\cdots f_r^{s_r} / \mathcal{D}_X[s_1,\dots,s_r]\cdot f_1^{s_1+1}\cdots f_r^{s_r+1},
$$
which for $r=1$ is the ideal generated by the Bernstein-Sato polynomial. The zero loci $V(B_F) \subset \mathbb{C}^r$ of Bernstein-Sato ideals are related with mixed multiplier ideals (see \cite{Mon21}) and cohomology support loci of local systems (see \cite{Bud15,BLSY17,BVWZ21}).

Let $X_1,\dots,X_r$ be smooth complex algebraic varieties, and $f_i \in \mathcal{O}_{X_i}(X_i)$ nonzero regular functions. Consider the product morphism 
$$
F=(f_1,\dots,f_r): X=\prod_{i=1}^r X_i \to \mathbb{A}^r.
$$ 
By \cite[Proposition 4.1]{BM99}, it follows that the Bernstein-Sato ideal $B_F$ is a principal ideal generated by the product of individual Bernstein-Sato polynomials
$$
B_F = (b_{f_1}(s_1)\cdots b_{f_r}(s_r)).
$$
Lemma \ref{main lemma} provides another proof of this result by taking $k:=\mathbb{C},A_i:=\mathbb{C}[s_i],P_i:=\mathcal{D}_{X_i}[s_i]f_i^{s_i},Q_i:=\mathcal{D}_{X_i}[s_i]f_i^{s_i+1},\delta_i:=f^{s_i},b_i:=b_{f_i}(s_i)$.

\subsection{Milnor fiber monodromy}

Let $X$ be a smooth complex algebraic variety and let $f\in \mathcal{O}_X(X)$ be a regular function. For $P \in X$ a point, the local Bernstein-Sato polynomial of $f$ at $P$ is defined as $b_{f,P}(s) := \gcd_U b_{f|_U}(s)$, where $U$ ranges over all open neighborhoods of $P$. The local Bernstein-Sato polynomials recover the global Bernstein-Sato polynomial via $b_f(s)=\lcm_{P \in X} b_{f,P}(s)$. By \cite[Corollary 3.34]{MJN21}, the Bernstein-Sato polynomial of $f$ and its analytification $f^{an} \in \mathcal{O}_{X^{an}}(X^{an})$ are equal. For a closed point $P \in V(f)$ in the vanishing locus of $f$, let $F_{f,P}$ denote the Milnor fiber of the map germ $f^{an}: (X^{an},P) \to (\mathbb{C},0)$ (see \cite{Mil68} for Milnor fibers), and denote the spectrum of eigenvalues of the monodromy action on the cohomology of $F_{f,P}$ by
$$
\Sp(f,P) := \bigcup_{i \in \mathbb{Z}} \{ \lambda \in \mathbb{C} \ | \ \text{$\lambda$ is an eigenvalue of the monodromy action on $H^i(F_{f,P}, \mathbb{C})$} \}.
$$
For $P \in X$ a closed point, we have by \cite[Theorem 6.3.5]{Bjo93} and its proof that
\begin{equation}\label{rel}
\bigcup_{P' \in U \cap V(f)} \Sp(f,P') = e^{2\pi i b_{f,P}^{-1}(0)}
\end{equation}
for $U$ any sufficiently small open neighborhood of $P$ in $X^{an}$ (also see \cite[Remark 2.12]{Saito_microlocal}, \cite[Remark 3.2]{Mas17}), where $b_{f,P}^{-1}(0)$ is the set of roots of $b_{f,P}(s)$.


Now let $X,Y$ be smooth complex algebraic varieties, and let $f\in \mathcal{O}_X(X), g \in \mathcal{O}_Y(Y)$ be regular functions, and consider the product function $fg \in \mathcal{O}_{X\times Y}(X \times Y)$. Fix closed points $P \in X, Q \in Y$. By Theorem \ref{main thm} we have the set-theoretic equality
\begin{equation}\label{roots}
b_{fg,(P,Q)}^{-1}(0) = b_{f,P}^{-1}(0) \cup b_{g,Q}^{-1}(0).
\end{equation}
In fact, there is an easier proof of this equality: it is immediate that $b_{fg}(s) \ | \ b_f(s)b_g(s)$, and 
$$
b_f(s)= b_{fg|_{X \times (Y\setminus V(g))}}(s) \ | \ b_{fg}(s),
$$
where the equality holds because $fg|_{X \times (Y\setminus V(g))}$ is the pullback of $f$ along a smooth morphism multiplied by an invertible function, and it is a standard fact that both these operations don't affect the Bernstein-Sato polynomial. 

By (\ref{rel}) and (\ref{roots}), we have 
\begin{equation*}
\bigcup_{(P',Q') \in (U\times W) \cap V(fg)} \Sp(fg,(P',Q')) = \bigcup_{P' \in U \cap V(f)} \Sp(f,P')  \cup \bigcup_{Q' \in W \cap V(g)} \Sp(g,Q'),
\end{equation*}
where $U \subset X^{an}$ and $W \subset Y^{an}$ are sufficiently small open neighborhoods of $P$ and $Q$. Here we must take the combined collection of Milnor monodromy eigenvalues over all points in sufficiently small open neighborhoods, so that we unfortunately don't get an explicit description of $\Sp(fg,(P,Q))$ in terms of $\Sp(f,P)$ and $\Sp(g,Q)$. However, by using the geometry of Milnor fibers in a special case, Tapp obtained the following pointwise relationship of Milnor monodromy eigenvalues (see \cite[Corollary 3.3.2]{Tap07} and \cite[Corollary 3.2]{Tap08}): if $f \in \mathbb{C}[x_1,\dots,x_n]$ and $g\in \mathbb{C}[y_1,\dots,y_m]$ are reduced homogeneous polynomials of possibly different positive degrees, then
$$
\Sp(fg,(0,0)) =  \Sp(f,0) \cap \Sp(g,0).
$$





\subsection{$V$-filtration}\label{vfilt}

We use the notation from \S\ref{setup}. The left $D_X[s]$-action on $B_{f,+}$ upgrades to a left action by the ring $D_{X \times \mathbb{A}^1}=D_X\langle t,\partial_t \rangle$ of differential operators on $X \times \mathbb{A}^1$, where $t$ and $\partial_t$ act by
$$
t \cdot us^jf^s = u (s+1)^j f^{s+1}, \ \ \ \partial_t \cdot us^jf^s =-u s(s-1)^j f^{s-1},
$$
where $u\in R_f$, and $j \in \mathbb{Z}_{\geq 0}$, so that $s=-\partial_t t$. If we let $i_f: X \to X \times \mathbb{A}^1, x\mapsto (x,f(x))$ be the graph embedding, then we get an isomorphism of left $D_{X \times \mathbb{A}^1}$-modules
$$
B_{f,+} \simeq i_{f,+}R_f = R_f[\partial_t]\delta_f,
$$
where $i_{f,+}$ is the $\mathcal{D}$-module theoretic direct image and the isomorphism is determined by sending $f^s$ to the formal symbol $\delta_f$. Define
$$
B_f := i_{f,+}R = R[\partial_t]\delta_f \subset B_{f,+}.
$$
Note that $D_{X\times \mathbb{A}^1}f^s=B_f \subset B_{f,+}$.

We now consider $V$-filtrations of $D_{X\times \mathbb{A}^1}$-modules along the smooth hypersurface $V(t)=X\times 0$ of $X \times \mathbb{A}^1$. For more details and properties of $V$-filtrations, see \cite{Sab87a}, \cite[\S 3.1]{Saito_MH},\cite[\S 2.2]{CD21}. The ring $D_{X \times \mathbb{A}^1}$ has a $\mathbb{Z}$-indexed $V$-filtration $(V^pD_{X \times \mathbb{A}^1})_{p\in \mathbb{Z}}$ given by
$$
V^pD_{X \times \mathbb{A}^1} = \sum_{i -j\geq p} D_Xt^i\partial_t^j.
$$
In particular, $V^0D_{X\times \mathbb{A}^1}=D_X\langle t,s\rangle$, 
$$
V^mD_{X\times \mathbb{A}^1}=t^mV^0D_{X\times \mathbb{A}^1}=V^0D_{X\times \mathbb{A}^1}t^m \ \ \ \ \text{for $m\in \mathbb{Z}_{\geq 0}$}
$$
and each $V^pD_{X\times \mathbb{A}^1}$ is a finitely generated left and right $V^0D_{X\times \mathbb{A}^1}$-module.

The $D_{X\times \mathbb{A}^1}$-module $B_{f,+}$ has a $\mathbb{Q}$-indexed $V$-filtration $(V^{\alpha}B_{f,+})_{\alpha \in \mathbb{Q}}$, which is the unique exhaustive, decreasing, and \textit{discrete and left-continuous}\footnote{This means that there is a positive integer $l$ such that $V^{\alpha}B_{f,+}$ has constant value for all $\alpha$ in an interval $(\frac{i-1}{l},\frac{i}{l}]$ with $i \in \mathbb{Z}$.} filtration satisfying the following properties:
\begin{enumerate}[(i)]
\item $V^pD_{X\times \mathbb{A}^1}\cdot V^{\alpha}B_{f,+} \subset V^{p+\alpha}B_{f,+}$ for every $p\in \mathbb{Z}, \alpha \in \mathbb{Q}$,
\item each $V^{\alpha}B_{f,+}$ is a coherent $V^0D_{X\times \mathbb{A}^1}$-submodule of $B_{f,+}$,
\item for every $\alpha \in \mathbb{Q}$, $s+\alpha$ acts nilpotently on $\Gr_V^{\alpha}B_{f,+}:= V^{\alpha}B_{f,+}/V^{>\alpha}B_{f,+}$, where $V^{>\alpha}B_{f,+}:=\cup_{\beta>\alpha} V^{\beta} B_{f,+}$, and
\item for every $p \in \mathbb{Z}_{\geq 0}$, we have $V^pD_{X\times \mathbb{A}^1}\cdot V^{\alpha}B_{f,+} = V^{p+\alpha}B_{f,+}$ for $\alpha \gg 0$.
\end{enumerate}
Similarly we have the $V$-filtration for the $D_{X\times \mathbb{A}^1}$-module $B_f$, which is given by $V^{\alpha}B_f = B_f \cap V^{\alpha}B_{f,+}$. Note that $V^{\alpha}B_f =  V^{\alpha}B_{f,+}$ for $\alpha>0$ by \cite[Lemma 3.1.7]{Saito_MH} because $B_{f,+}=B_f[1/t]$. 

Similarly consider the $V$-filtrations on $B_{g,+}$ and $B_{fg,+}$. 

Recall from (\ref{tensor}) the following isomorphism over $D_{X\times Y}[s] = D_X[s] \otimes_{\mathbb{C}[s]} D_Y[s]$:
$$
B_{fg,+}\simeq B_{f,+}\otimes_{\mathbb{C}[s]} B_{g,+}.
$$
The left hand side is a $D_{X \times Y\times \mathbb{A}^1}$-module and the right hand side is a $D_{X \times \mathbb{A}^1} \otimes_{\mathbb{C}[s]}D_{Y \times \mathbb{A}^1}$-module, so we indicate below how the $t$ and $\partial_t$ actions correspond on both sides:
\begin{align*}
D_{X \times Y\times \mathbb{A}^1}[t^{-1}] \supset D_{X\times Y}[s] &\simeq D_X[s] \otimes_{\mathbb{C}[s]} D_Y[s] \subset D_{X \times \mathbb{A}^1}[t^{-1}] \otimes_{\mathbb{C}[s]}D_{Y \times \mathbb{A}^1}[t^{-1}] \\
s &\leftrightarrow s \otimes 1 = 1 \otimes s \\
t &\leftrightarrow t \otimes t \\
\partial_t = -st^{-1} &\leftrightarrow -(s \otimes 1)(t^{-1} \otimes t^{-1}) = -st^{-1} \otimes t^{-1} = \partial_t \otimes t^{-1}.
\end{align*}


\begin{lemma}\label{vfilt tensor}
Under the isomorphism $B_{fg,+} \simeq B_{f,+}\otimes_{\mathbb{C}[s]} B_{g,+}$ we have 
$$
V^{\alpha}B_{fg,+} \simeq V^{\alpha}B_{f,+} \otimes_{\mathbb{C}[s]} V^{\alpha}B_{g,+}  \ \ \ \text{for all $\alpha \in \mathbb{Q}$}.
$$
\end{lemma}
\noindent{\it Proof}. Let $W^{\alpha} = V^{\alpha}B_{f,+} \otimes_{\mathbb{C}[s]} V^{\alpha}B_{g,+}$. By uniqueness of the $V$-filtration, it suffices to show that $W^{\bullet}$ satisfies the properties of the $V$-filtration for the $D_{X \times Y \times \mathbb{A}^1}$-module $B_{fg,+}$. It is easy to verify that $W^{\bullet}$ satisfies the properties of the $V$-filtration, except perhaps the finite generation over $V^0D_{X \times Y \times \mathbb{A}^1}$, which we show as follows. Observe that
$$
V^0D_{X\times \mathbb{A}^1} f^s = D_X[s]f^s , V^0D_{Y\times \mathbb{A}^1} g^s = D_Y[s]g^s, V^0D_{X\times Y\times \mathbb{A}^1} (fg)^s = D_{X\times Y}[s](fg)^s,
$$
so that
$$
V^0D_{X\times Y\times \mathbb{A}^1}(fg)^s \simeq V^0D_{X\times \mathbb{A}^1}f^s \otimes_{\mathbb{C}[s]} V^0D_{Y\times \mathbb{A}^1}g^s.
$$
Fix $\alpha \in \mathbb{Q}$. Note that if $W^{\alpha}$ is generated over $V^0D_{X\times Y\times \mathbb{A}^1}$ by some elements $u_i$, then $tW^{\alpha}$ is generated by $tu_i$ over $V^0D_{X\times Y\times \mathbb{A}^1}$ because
$$
tW^{\alpha} = t\sum_i V^0D_{X\times Y\times \mathbb{A}^1}u_i =\sum_i V^1D_{X\times Y\times \mathbb{A}^1}u_i=  \sum_i V^0D_{X\times Y\times \mathbb{A}^1}tu_i
$$
and the converse is also true because the $t$-action is invertible on $B_{fg,+}$. So it suffices to show that  $t^mW^{\alpha}$ is a coherent $V^0D_{X\times Y\times \mathbb{A}^1}$-module for $m \gg 0$. Because $t^mW^{\alpha} \subset W^{\alpha+m}$ for $m\geq 0$ and $V^0D_{X\times Y\times \mathbb{A}^1}$ is a noetherian ring, it suffices to show that $W^{\alpha+m}$ is a coherent $V^0D_{X\times Y\times \mathbb{A}^1}$-module, hence we may assume $\alpha > 0$. 

Because $B_f = D_{X\times \mathbb{A}^1}f^s = \cup_m V^mD_{X\times \mathbb{A}^1}f^s$ is an exhaustive increasing union by coherent $V^0D_{X\times \mathbb{A}^1}$-modules, by coherence of $V^{\alpha}B_{f}$ over $V^0D_{X\times \mathbb{A}^1}$ we have
$$
V^{\alpha}B_{f,+}=V^{\alpha}B_{f} \subset V^{-m}D_{X\times \mathbb{A}^1}f^s \ \ \ \ \text{for $m \gg0$},
$$ 
hence for $m\gg 0$ we have
$$
t^mV^{\alpha}B_{f,+}\subset t^mV^{-m}D_{X\times \mathbb{A}^1}f^s \subset V^0D_{X \times \mathbb{A}^1}f^s
$$
and
$$
t^mW^{\alpha} =  t^mV^{\alpha}B_{f,+} \otimes_{\mathbb{C}[s]} t^mV^{\alpha}B_{g,+} \subset V^0D_{X \times \mathbb{A}^1}f^s \otimes_{\mathbb{C}[s]}V^0D_{Y \times \mathbb{A}^1}g^s \simeq V^0D_{X\times Y\times \mathbb{A}^1}(fg)^s.
$$
Thus $t^mW^{\alpha} $ is finitely generated over $V^0D_{X \times Y \times \mathbb{A}^1}$ by noetherianity of $V^0D_{X \times Y \times \mathbb{A}^1}$. \qed

\subsection{Generalized $b$-functions}

For an element $u \in B_{f,+}$, the generalized $b$-function $b_u(s)\in \mathbb{C}[s]$ of $u$ is the minimal polynomial of the action of $s$ on $V^0D_{X \times \mathbb{A}^1}\cdot u/ V^1D_{X \times \mathbb{A}^1}\cdot u$. We recover the Bernstein-Sato polynomial of $f$ by $b_{f^s}(s)=b_f(s)$. By Sabbah's description of the $V$-filtration (see \cite{Sab87a}), we have
$$
V^{\alpha}B_{f,+} = \{ u \in B_{f,+} \ | \ \text{the roots of $b_u(-s)$ are $\geq \alpha$} \}.
$$
For $u\in B_{f,+}, v \in B_{g,+}$, and their product $u\otimes v \in B_{fg,+}$, we can ask whether there is a relationship between $b_{u\otimes v}(s)$ and $b_u(s)b_v(s)$ as we have found for $u=f^s, v=g^s$. By the same proof, we have
$$
b_{pf^s}(s)b_{qg^s}(s) = b_{pq(fg)^s}(s) \ \ \ \text{for all $p \in R_f, q\in S_g$}.
$$
But for general $u,v$, neither side needs to divide the other, as the example below shows. The most we can say in general is that because $V^{\alpha}B_{fg,+} \simeq V^{\alpha}B_{f,+} \otimes_{\mathbb{C}[s]} V^{\alpha}B_{g,+}$, by Sabbah's description of the $V$-filtration, if all the roots of $b_u(-s)$ and $b_v(-s)$ are $\geq \alpha$, then all the roots of $b_{u \otimes v}(-s)$ are also $\geq \alpha$ (in fact, this is equivalent to saying $V^{\alpha}B_{f} \otimes_{\mathbb{C}[s]} V^{\alpha}B_{g}$ is contained in $V^{\alpha}B_{fg}$).



Before we give the example where neither $b_{u\otimes v}(s)$ nor $b_u(s)b_v(s)$ divides the other, we make a few observations.
\begin{enumerate}[(i)]
\item Let $u\in B_{f,+}$. Then $b_{t^mu}(s)=b_u(s+m)$ for all $m \in \mathbb{Z}$ because the action of $t$ is invertible on $B_{f,+}$.
\item For this observation, we follow the same argument as in \cite[Proposition 7.12]{MP20}. Let $u\in R_f$, and assume we have the functional equation
$$
(s+1)b(s) uf^s = P(s)tuf^s \in D_X[s]tuf^s
$$
for some $b(s) \in \mathbb{C}[s]$ and $P(s) \in D_X[s]$. Applying $t^{-1}$ to both sides, we get 
$$
-\partial_tb(s)uf^s = P(s-1)uf^s.
$$
Now observe that
$$
(s+1)b(s-m)u\partial_t^mf^s = (s+1)u\partial_t^mb(s)f^s = -t\partial_tu\partial_t^mb(s)f^s
$$
$$
=  -tu\partial_t^m\partial_tb(s)f^s = tu\partial_t^mP(s-1)uf^s = P(s-m) tu\partial_t^mf^s \in D_X[s] t u\partial_t^mf^s.
$$
The takeaway is
$$
b_{u\partial_t^mf^s}(s) \ | \ (s+1)b(s-m).
$$
\item If $f$ is not an invertible element of $R$, then $\cap_{\alpha} V^{\alpha}B_{f,+}=0$. To see why, first note that $\cap_{\alpha} V^{\alpha}B_{f,+} \subset B_f$ because $V^{\alpha}B_{f,+} = V^{\alpha}B_f$ for $\alpha>0$. Since $t$ is an isomorphism on $B_{f,+}$, it follows that $t$ is an isomorphism on $\cap_{\alpha} V^{\alpha}B_{f,+}$, hence the intersection is a module over $V^0D_{X\times \mathbb{A}^1}[t^{-1}]=D_{X\times \mathbb{A}^1}[t^{-1}]$. But by Kashiwara's equivalence, $B_f=i_{f,+}R$ is a simple $D_{X\times \mathbb{A}^1}$-module, hence $\cap_{\alpha} V^{\alpha}B_{f,+}$ is either all of $i_{f,+}R$ or zero, hence by Sabbah's description of the $V$-filtration, $\cap_{\alpha} V^{\alpha}B_{f,+}$ is zero for $f$ not a unit in $R$.
\item The jumps of the $V$-filtration on $B_{f,+}$ coincide precisely with the roots of $b_f(-s)$ modulo $\mathbb{Z}$. In particular, if all the roots of $b_f(s)$ are integers, then $V^{\alpha}B_{f,+}=V^{\lceil \alpha \rceil}B_{f,+}$ for all $\alpha \in \mathbb{Q}$.
\end{enumerate}

\begin{example}
Let $R=\mathbb{C}[x],S=\mathbb{C}[y]$ be one variable polynomial rings, $f=x,g=y$, and $u=x\partial_t^2x^s, v= t^{-2}y^s$, and consider the product
$$
u\otimes v = x\partial_t^2 (xy)^s  \in B_{fg,+} \simeq B_{f,+}\otimes_{\mathbb{C}[s]} B_{g,+}.
$$
We claim that neither $b_{x\partial_t^2 (xy)^s}$ nor $b_{x\partial_t^2x^s}(s)b_{t^{-2}y^s}(s)$ divides the other. For this, it suffices to show
\begin{align}
b_{t^{-2}y^s}(s) = s-1 \\
b_{x\partial_t^2x^s}(s) =s+1 \\
s \  | \  b_{x\partial_t^2 (xy)^s}(s) \ | \ (s+1)s. \label{moyp}
\end{align}
The first equation follows from observation (i). The second equation is shown as follows. Applying $(s+1)t$ to $(s+1)x^s = \partial_x x^{s+1}$, we get $(s+1)(s+2)xx^s \in D_X[s]txx^s$, hence $b_{x\partial^2_t x^s}(s) \ |\ (s+1)s$ by observation (ii), but $x\partial_t^2x^s=x\partial^2_x x^s \in V^1$ (because $\partial_t x^s=-\partial_x x^s$), so by Sabbah's description of the $V$-filtration, we must have $b_{x\partial^2_t x^s}(s)=s+1$ (it cannot be $1$ by observation (iii)).

For (\ref{moyp}), first note that the multiplicative Thom-Sebastiani property applies to $x(xy)^s=xx^s \otimes y^s$ so that
$$
b_{x(xy)^s}(s) = b_{tx^s \otimes y^s}(s)=b_{tx^s}(s)b_{y^s}(s) = (s+2)(s+1).
$$
Because $V^1D_{X\times Y \times \mathbb{A}^1}\cdot x(xy)^s = D_{X\times Y}[s] xt(xy)^s$, we can use observation (ii) to say
$$
b_{x\partial_t^j(xy)^s}(s) \ |\  (s+1)(s+2-j) \ \ \ \text{for all $j\geq 0$}.
$$
To show that $s$ divides $b_{x\partial_t^2(xy)^s}(s)$, it suffices to show $b_{x\partial_t^2(xy)^s}(s)$ is not equal to $1$ or $s+1$. Since the $V$-filtration is separated for noninvertible functions, $b_{x\partial_t^2(xy)^s}(s)\neq 1$. To show $s$ divides $b_{x\partial_t^2(xy)^s}(s)$, it suffices to show $x\partial_t^2(xy)^s\in V^0B_{xy,+} \setminus V^1B_{xy,+}$. By Sabbah's description of the $V$-filtration, $x\partial_t(xy)^s\in V^1B_{xy,+}$ and $x\partial_t^2(xy)^s\in V^0B_{xy,+}$. Observe that
$$
V^1B_{xy,+} \ni \partial_x x\partial_t(xy)^s = \partial_t(xy)^s - xy\partial_t^2(xy)^s.
$$
Note that $\partial_t(xy)^s \not\in V^1B_{xy,+}$ because $s | b_{\partial_t(xy)^s}(s)$ by  \cite[Proposition 6.11]{MP20}. Thus we must also have $xy\partial_t^2(xy)^s \not\in V^1B_{xy,+}$, hence $x\partial_t^2(xy)^s\not\in V^1B_{xy,+}$ and we are done.
\end{example}








\end{document}